\documentclass{amsart}

\usepackage{amsmath}
\usepackage{amsthm}
\usepackage{amssymb}
\usepackage[all]{xy}
\usepackage[mathscr]{eucal}

\UseComputerModernTips

\usepackage{amsthm}

\newcommand{\rat}{\ensuremath{\mathbb{Q}}}
\newcommand{\ganz}{\ensuremath{\mathbb{Z}}}
\newcommand{\tate}{\ensuremath{\mathbb{L}}}
\newcommand{\Knull}{\ensuremath{\operatorname{K}_0}}

\newcommand{\Vark}{\ensuremath{\mathit{Var}_k}}
\newcommand{\reg}{\ensuremath{\mathcal{O}}}
\newcommand{\Fq}{\ensuremath{\mathbb{F}_q}}
\newcommand{\aff}[1]{\ensuremath{\mathbb{A}^{#1}}}
\newcommand{\proj}[1]{\ensuremath{\mathbb{P}^{#1}}}
\newcommand{\Spec}{\ensuremath{\operatorname{Spec}}}

\newcommand{\SBk}{\ensuremath{\operatorname{SB}_k}}

\newtheorem{thm}{Theorem}[section]

\newtheorem{qst}[thm]{Question}

\begin{document}
\title{A note on congruences for theta divisors}
\author{Franziska Heinloth}
\address{Universit\"at Duisburg---Essen, Standort Essen, 
FB6, Mathematik, 45117 Essen, Germany}
\email{franziska.heinloth@uni-duisburg-essen.de}
\begin{abstract} The classes of two theta divisors on an abelian variety in the naive Grothendieck ring of varieties need not be congruent modulo the class of the affine line.
\end{abstract}
\maketitle

\section{Introduction}
A \emph{theta divisor} on an abelian variety $A$ over a field is an effective ample divisor $\Theta$ such that the global sections of the line bundle $\reg_A(\Theta)$ are one-dimensional. For example, the image of the $(g-1)$st symmetric power of a smooth projective curve $C$ of genus $g$ with a rational point in its Jacobian is a theta divisor.

Berthelot, Bloch and Esnault prove the following 
(\cite{BerthelotBlochEsnault}, Theorem 1.4):

\begin{thm}
Let $\Theta$, $\Theta'$ be two theta divisors defined over a finite field $\Fq$.
Then $$|\Theta(\Fq)|\equiv|\Theta'(\Fq)| \mod q.$$
\end{thm}

This was conjectured by Serre, who more generally predicted that in 
some sense to be made precise, on an abelian variety defined over a field $k$, the motives of two theta divisors $\Theta$ and $\Theta'$ should differ by a multiple of the Lefschetz motive.

A naive version of this conjecture is the following

\begin{qst}
Is $[\Theta]\equiv[\Theta'] \mod \tate$ in $\Knull(\Vark)$?
\end{qst}  

Here $\Knull(\Vark)$ denotes the naive Grothendieck ring of varieties over $k$, i.e. the free abelian group on isomorphism classes $[X]$ of varieties over $k$ modulo the relations $[X]=[X-Y]+[Y]$ for $Y\subset X$ a closed subvariety. The ring structure is given by the product of varieties, and $\tate$ is the class $[\aff{1}_k]$ of the affine line.

In this note we give a negative answer to the above question using the fact that products of two elliptic curves can be Jacobians.

\section{A counterexample on a product of two elliptic curves}
In \cite{HayashidaNishi}, Hayashida and Nishi show that certain products of two elliptic curves $E\times E'$ are Jacobians. More precisely they prove (Theorem 4 and Proposition 4):
\begin{thm}
\label{kurve}
Let $E$ and $E'$ be isogenous elliptic curves with endomorphism ring $\reg$ the ring of integers in the number field $\rat(\sqrt{-m})$, where $m$ is a nonnegative integer.
Then $E\times E'$ is a Jacobian variety of a curve of genus $2$ if $E'$ is not isomorphic to $E$ or if $m \not\in\{0,1,3,7,15\}$.
\end{thm}

Recall that two smooth projective varieties $X$ and $X'$ are called \emph{stably birational}, if $X\times\proj{n}$ and $X'\times\proj{n'}$ are birational for some integers $n$ and $n'$. Denote the set of stably birational classes of irreducible smooth projective varieties over $k$ by $\SBk$, denote the free abelian group on $\SBk$ by $\ganz[\SBk]$. It carries a commutative ring structure induced by the product of varieties.

The weak factorization theorem by W{\l}odarczyk and Abramovich et al. (see \cite{Wlodarczyk} and \cite{Abramovich}) implies that two birational irreducible smooth projective varieties over a field of characteristic zero are connected by a series of blow-ups and blow-downs along smooth centers.
Larsen and Lunts in \cite{LarsenLunts}, using the weak factorization theorem, show that over a field $k$ of characteristic zero the following holds:

\begin{thm}
There is a ring homomorphism $b:\Knull(\Vark)\longrightarrow \ganz[\SBk]$ sending the class of an irreducible smooth projective variety $[X]$ to its stably birational class $[X]_b$.
\end{thm}

Obviously, this homomorphism sends the class $\tate$ of the affine line to zero. On the other hand, using the weak factorization theorem again, it is easy to see that two stably birational irreducible smooth projective varieties have the same class in $\Knull(\Vark)/(\tate)$. Therefore, as pointed out in \cite{LarsenLunts}, 
$b$ actually induces an isomorphism of $\Knull(\Vark)/(\tate)$ and $\ganz[\SBk]$.

Now take $\Theta$ a genus two curve with a rational point over a field $k$ of characteristic zero such that the Jacobian of $\Theta$ is $E\times E'$ as in Theorem \ref{kurve}. Let $\Theta'= E\times {0'} + 0\times {E'}$. Then
$b([\Theta])=[\Theta]_b\ne b([\Theta'])=[E]_b+[E']_b-[\Spec k]_b$, as a curve of genus $2$ is not stably birational to an elliptic curve (note that we are working in a \emph{free} abelian group).

Therefore the classes of $\Theta$ and $\Theta'$ are not congruent modulo $\tate$ in the Grothendieck ring of varieties.

\subsection*{Acknowledgments} I thank H\'el\`ene Esnault for drawing my
attention to congruence questions for theta divisors and for helpful comments.

\end{document}